\documentclass[journal]{IEEEtran}
\usepackage{amsmath}
\usepackage{amssymb}
\usepackage{amsthm}
\usepackage{graphicx}
\usepackage{textcomp}
\usepackage{cite}
\usepackage{url}
\usepackage{array}
\usepackage{color,soul}
\usepackage{bm}
\usepackage{epstopdf}

\begin{document}

\title{{\bf S}mart {\bf P}rocurement {\bf O}f {\bf N}aturally {\bf G}enerated {\bf E}nergy (SPONGE) for PHEV's \thanks{This work was supported in part by Science Foundation Ireland grant 11/PI/1177.}}

\author{Florian~H\"{a}usler, Emanuele~Crisostomi, Ilja~Radusch and Robert~Shorten%
\thanks{F. H\"{a}usler and I. Radusch are with the Fraunhofer Institute for Open Communication Systems, 10623 Berlin, Germany.}%
\thanks{E. Crisostomi is with the Department of Energy, Systems, Territory and Constructions Engineering, University of Pisa, L.go L. Lazzarino, 56122, Pisa, Italy.}%
\thanks{R. Shorten is with IBM Research Ireland, Dublin 4, Ireland and Hamilton Institute, NUI Maynooth, Ireland}%
}

\maketitle

\begin{abstract}
In this paper we propose a new engine management system for hybrid vehicles to enable energy providers and car manufacturers to provide new services. Energy forecasts are used to collaboratively orchestrate the behaviour of engine management systems of a fleet of PHEV's to absorb oncoming energy in an smart manner. Cooperative algorithms are suggested to manage the energy absorption in an optimal manner for a fleet of vehicles, and  the mobility simulator SUMO is used to show simple simulations to support the efficacy of the proposed idea.

\end{abstract}

\section{Introduction}
\label{Introduction}
One of the promised advantages of the Smart Grid is the ability to integrate power generated from renewable sources into the daily demands of society in a manner that takes into account the intrinsic fluctuating nature of such intermittently available power. In spite of the challenges in dealing with an uncertain supply, several countries have already started producing a large fraction of electrical power from renewable sources. For example, wind alone provided more than 30~\% of electricity production in Denmark in 2012, and is foreseen to supply 50~\% of the overall demand by the year 2020~\cite{Marinelli}. Further, Denmark's stated goal goes beyond this objective with an aspiration to become 100~\% renewable by 2050~\cite{Meibom}, with similar aspirations being held in many other countries. The merits of using renewables is the cleanness of the energy supply and the associated benefits for both greenhouse gas emissions and air quality. \newline 

A significant impediment to the integration of renewables into the grid is the need for new demand side management practices to match power generation with power consumption~\cite{Bicik}, \cite{RAMP} on a daily basis. Despite the increasing quantity of energy that is produced from renewable sources, and despite the many efforts to encourage consumers to shift loads to times of the day when renewable energy is available, there is still a significant mismatch between renewable energy availability and energy demand, and conventional power plants (e.g., coal or gas-fired power plants) are still widely used to back-up energy generation. The necessity to use conventional power plants is not convenient in terms of economic costs (i.e., fuel and carbon costs have to be taken into account), and in terms of existing and anticipated environmental regulations (e.g., emissions of $CO_2$, NOx or other pollutants).\newline

With this in mind, a number of strategies have been proposed to deal with supply-demand imbalance.  First, storage systems represent an attractive possibility to alleviate the requirement of continuous matching between energy demand and offer; see the Economist Technology Quarterly (December 2014) for a recent discussion of advances in this direction \cite{Economist}.  Roughly speaking, energy generated from renewable sources can be stored when availability exceeds the energy demand (e.g., eolic energy at night time), and can be released as needed as an alternative to switching on a conventional power plant. Amongst the available storage systems, the ability of Electric Vehicles (EVs) to act as an 'aggregated' battery for such purposes has been given as one of the most important arguments in favour of EV adoption as a mode of automotive mobility.  To further elaborate: in the event of a high level of adoption of electric vehicles, a large virtual battery system would be automatically available without requirements of big investment in other storage devices ~\cite{Liu},~\cite{Tushar} and \cite{Griggs}. Such a possibility would be convenient also from the EV side, as many studies show that EVs represent a viable solution to limit emissions particularly if they are charged from energy coming from renewable sources \cite{PNAS}. Despite the apparent suitability of EVs to provide battery storage, some studies have shown that the use of EVs in this manner is not without problems. Their use as storage devices is currently quite limited due to low penetration levels, and due to the fact that, usually, only a small proportion of the EVs' battery capacity is available for energy exchange \cite{CNR_IEVC} (the main reason for this is that EV owners are primarily concerned on having enough battery for their next trip \cite{Griggs}).  A further complication arises due to the fact that renewable energy is, by its very nature, uncertain. Thus, supplying the EV fleet with the requited energy needed for mobility creates the need for complicated optimisation and scheduling algorithms (and infrastructure) on the supply side,  places contractual requirements on generators of electricity, and requires EV owners to plug in their vehicles at certain times of the day \cite{Edash}.\newline
 
A second strategy in dealing with the intermittent nature of the supply of naturally generated energy is to make devices smarter. That is, they should be context aware, and change their behaviours in a manner that enables them to utilise renewable energy when it becomes available. Clearly, in the case of EVs this is not possible since there is a very natural decoupling of the mobility needs of EV owners and the available supply of energy from renewable sources. That is, the grid should always serve the mobility needs of users as a primary objective. However, in the case of PHEVs, there are two power sources for each vehicle. Thus, vehicle owners have a choice at every instant of time as to whether the EV engines is utilised, the ICE, or both. Our objective in this paper is to describe how this flexibility can be exploited to couple the needs of the grid with that of the vehicle owner, and to show how this offers the potential for a truly smart integration of vehicles into the energy grid. Specifically, using energy forecasting models, it is possible for PHEVs to act as an energy sponge and be primed to capture renewable 
energy as it becomes available. To do this weather forecast services can be used to make predictions of how much energy will be available from solar/wind power plants in the near future (e.g., next 24 hours). Based on these forecasts, elementary cooperative strategies can be implemented, in a manner that is transparent to users, to make space in a fleet of vehicles (in the batteries) for forthcoming energy. Essentially, vehicle owners allow the EMU (engine management unit) of the vehicle to be orchestrated by a management service that takes into account the needs of the grid. By doing this, users fully utilise the available clean energy as it becomes available (possibly at zero financial cost if coming from private roof top solar panels or from self-owned wind source),  prevent clean energy from being wasted by ensuring that there is always enough capacity available for storage, and help balance energy supply and demand through active scheduling of energy sourcing for vehicles in a pro-active manner. We shall also see that this strategy has the potential to significantly reduce the complexity burden of charging these vehicles by enabling best effort charging algorithms to be deployed.\newline 

This work presented here complements previous work by the authors in the context of optimised charging of EV's and in developing context cooperative control strategies for hybrid vehicles. From a technological perspective, the work most resembles strategies to regulate pollution that have already been implemented in practice and described \cite{Twinlin}. However, the goals of the current problem statement differ significantly from the aforementioned work, and open, we believe, significant market opportunities by presenting mobility and energy products in a truly integrated manner.

\section{Problem statement and Assumptions}
\label{Problem_Formulation}
For convenience, and ease of exposition, we make the following set of simplistic assumptions. We start by indexing each day (a 24 hour period) with the index $k$ (meaning the $k$'th day or time period). We divide the $k$'th day into two time periods; a time period when vehicles are charging, and a time period when vehicles are not charging and {\em possibly} in transit. Furthermore:\newline
\begin{itemize}
\item[(i)] We assume a group of $N \in (1,...,n_{max})$ plug-in hybrid vehicles that participate in a scheme to absorb available naturally generated energy. Alternatively, in special cases, if needed, they might also provide energy to the grid. Although the overall discussion goes along the same lines, for the sake of simplicity in what follows we shall primarily consider the G2V (Grid-to-Vehicle) case, although the V2G (Vehicle-to-Grid) case can be considered as well.\newline
\item[(ii)] Each vehicle is assumed to be capable of operating in fully EV mode, in ICE mode, or a combination of both (as in the Toyota Prius).\newline
\item[(iii)] We also assume that for some fixed period during the day, these vehicles are plugged in, and that for this period, a reliable day-ahead forecast of available renewable energy is available. We denote this available energy by $E_{av}(k)$. For example, a typical assumption might be that that the vehicles charge from 11pm to 6am, though it is not necessary for this time period to be the same for all vehicles. Although, in principle, the future horizon of optimisation can be longer than one day, weather forecasts might not be reliable enough to support optimal decisions over longer time periods, see \cite{Wind}-\cite{Sun}.\newline
\item[(iv)] For the remainder of the day, we assume that vehicles may be in service, at any instant of time $n(t) \leq N$ vehicles are in transit, and that these vehicles can report their energy consumption over some period to a central agent.\newline
\end{itemize}

\subsection{Smart Procurement of Energy: SPONGE}
\label{Sponge_case}
\begin{figure}[ht]
\centering
\includegraphics[width=\columnwidth]{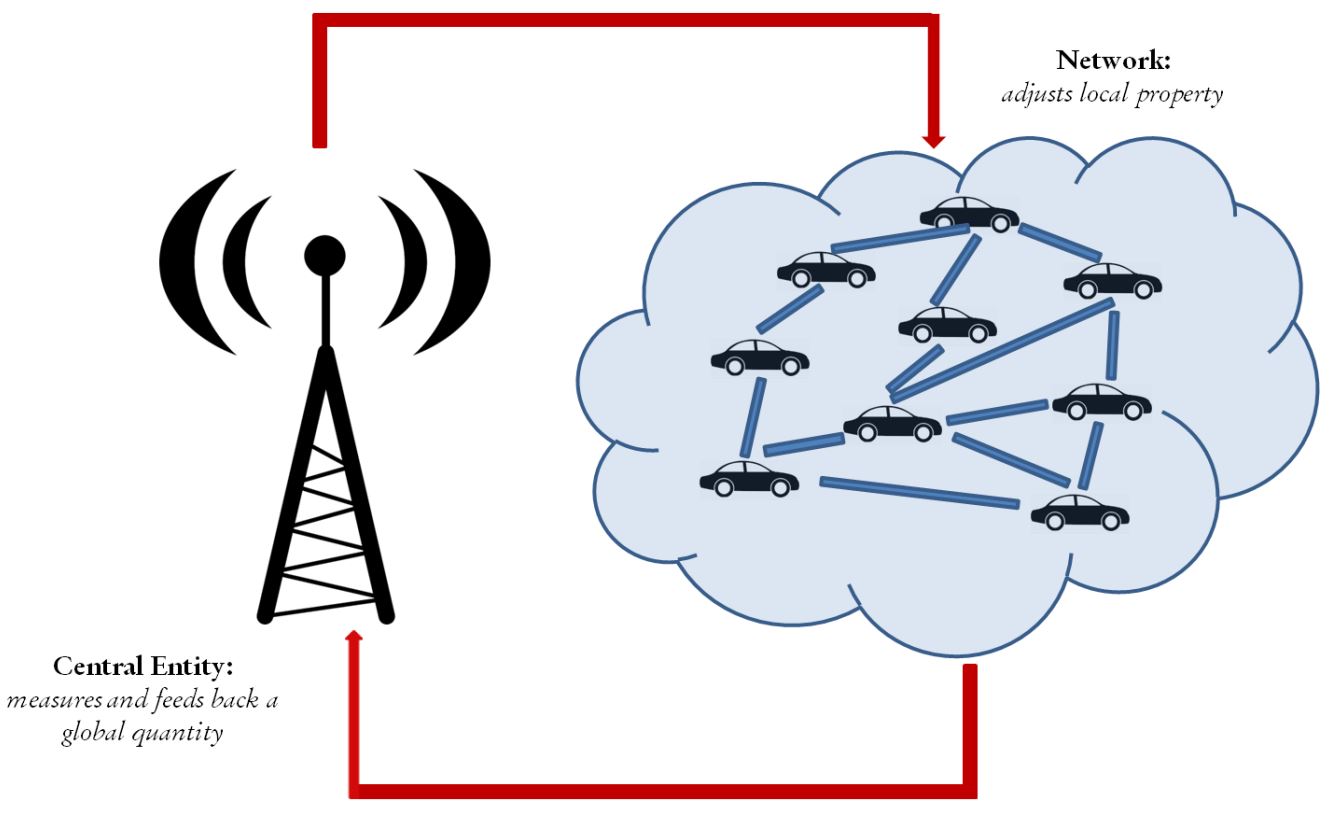}
\label{fig:1}
\caption{Feedback loop for energy dissipation problem}
\end{figure}

Let us now denote the electric energy dissipated by the $i$'th vehicle by $D_i(k)$. Our objective is to ensure that
\begin{eqnarray}
\sum_{i=1}^{N} D_i(k) \geq E_{av}(k+1).
\label{SPONGE}
\end{eqnarray} 
That is, during the $k$'th day the fleet acts like a sponge and makes available {\em at least} enough space to absorb the available energy that is expected 
during the next charging period. As stated, the problem is essentially a regulation problem that is depicted in Figure \ref{fig:1}. Under ideal circumstances,
a central authority computes the desired electrical energy consumption, and then broadcasts some signal which is received by the vehicles to orchestrate the switching between EV and ICE mode, so as to satisfy the regulation constraint. For instance, the signal can be the probability to travel in EV mode rather than in ICE mode, or can be the proportion of the traction torque that should be provided by the EV engine rather than from the ICE engine. We shall denote the problem expressed by Equation \ref{SPONGE} as the basic SPONGE problem.

\subsection{Smart Procurement of Energy: Exact SPONGE}
\label{Exact_Sponge_Case}
In some cases, the objective can be to make PHEVs travel in EV mode until they deplete their batteries in order to exactly match the expected energy that will be available from renewable sources. We shall denote this problem as ``exact SPONGE'', and its mathematical formulation is as follows:
\begin{eqnarray}
\sum_{i=1}^{N} D_i(k) = E_{av}(k+1).
\label{Exact_Sponge}
\end{eqnarray} 
The main advantage of the exact SPONGE approach is that when the fleet of vehicle connect to the grid for recharging, the quantity of required energy is already known in advance (i.e., it is equal to the expected energy available from renewable sources).

\subsection{Optimised access: Optimal SPONGE}
\label{Utility_Functions}
In some situations, certain vehicles may have prioritised access to the oncoming energy $E_{av}(k+1)$ via some utility function $f_i(D_i(k))$. Thus, the above problem can be reformulated in an optimisation framework as:
\begin{equation}\left\{
\begin{array}{lll}
\mbox{minimize} & & \displaystyle\sum_{i=1}^{N} f_i(D_i(k)) \\
& & \\
\mbox{subject to}   & & \displaystyle\sum_{i=1}^{N} D_i(k) = E_{av}(k+1).\end{array}\right.
\label{Utility_Maximisation_Problem}
\end{equation} 
This optimisation may be solved in many ways under suitable assumptions on the $f_i(D_i(k))$. The problem is most interesting when the the $f_i's$ represent a generalised notion of utility (or price that the $i$'th car pays) and is considered to be private information, not to be revealed to the utility or to other vehicles. The problem is then to solve the problem in a privacy preserving manner. Note that the $f_i's$ may be incorporated to represent various use cases. Some interesting examples include the following.
\begin{itemize}
\item[(i)] For example, OEM's may partner with utilities to provide a service where the price of energy is part of PHEV's owners car purchase plans. Those paying more upfront, may have prioritised access to 'free energy' as it becomes available \cite{Edash}.\newline 
\item[(ii)] The $f_i$'s could represent the price paid by an individual vehicle owner for energy access.\newline 
\item[(iii)] Or, they could be used to penalise vehicles with a lower load factor (fewer passengers).\newline
\item[(iv)] They could be used to penalise vehicles that drive close to schools, hospitals, etc.\newline 
\item[(v)] Another interesting scenario is as follows. Some hybrid modes blend the EV motor with the ICE to optimise fuel economy/emissions. An interesting embodiment of the optimisation scenario is to take the required energy in a manner that minimises the impact on fuel economy of the fleet.\newline
\end{itemize}

With regard to the SPONGE formulation several comments are appropriate.\newline

{\bf Comment 1: } Note that the SPONGE solution has the potential to completely revolutionise the ``charging paradigm''. Hitherto, most charging research has focussed on how to share the available energy among the connected fleet of vehicles in a manner that is compliant with the desires of the EV owners, the constraints of the grid, and the available power. Note that in this case, there might arise some problems in the power grid to accept the unexpected load, with the ultimate possibility of causing thermal overload of network components, low voltages at sensitive locations of the network, and increased phase unbalance \cite{deHoog}. Even ignoring this, the required optimisations often place severe constraints on the EV owners in the form of inconvenient charging profiles. On the other hand, in the solution of Problem (\ref{Exact_Sponge}), one would compute the same quantity in advance, and deplete the batteries of the vehicles while travelling of the same quantity. Thus, the charging process becomes fully schedulable and programmable. The charging problem can be reduced to a bast-effort problem where the cars share the available energy during the charging period using some simple algorithm such as Additive Increase Multiplicative Decrease (AIMD) algorithms \cite{Sonja,mingming}. Thus,  clearly, the difficulties of matching the demand and the offer are shifted to the driving stage through an optimal orchestration of the ICE and EV engines.\newline

{\bf Comment 2: } The discerning reader may ask why the individual vehicle owners should not simply expend the electric energy completely before switching to ICE mode. There are many reasons for doing this.  First, in some engines, electrical power and ICE are combined to reduce overall consumption, or for other objectives of interest (e.g., extend the lifetime of the battery) \cite{Oak_Ridge}.
Thus, it is advantageous to keep a store of naturally generated electrical energy for this purpose. Second, access to certain parts of the city may be restricted to zero emission vehicles. Thus, maintaining a store of electrical energy for this purpose is also advantageous. Finally, depleting the battery beyond the energy levels available during the next charging period, may lead to a situation where the battery is not filled during the $k+1$'th charging period. Thereby, the ICE may need to be engaged prematurely in driving, thus leading to unnecessary emissions and increased fuel consumption.\\
\newline
{\bf Comment 3: } Note that in some cases, depending on the number of vehicles on the road, the previous optimisation problems might not have a feasible solution. For instance, in the particular case that there are no vehicles on the road, then obviously the PHEVs can not deplete their batteries to make room for the forthcoming energy. In such cases where the problem does not have a solution, we will be interested in a `best-effort' solution, where the closest feasible solution is achieved instead, see for instance \cite{Sonja}.

\section{Methods}
\label{optimisation_methods}
Clearly there are many ways in which the problems specified in the previous section may be solved, and we now describe some simple methods that can be adopted. To this end we assume that over the period when vehicles are not charging, denoted $\theta(k)$, the aggregate electrical energy consumption of the fleet is required to be $E(k,t)$ with, 
\begin{eqnarray}
\int_{\theta(k)} E(k,t) dt = E_{av}(k+1)
\end{eqnarray}
where $t$ denotes time. We also assume that each vehicle is synchronised with a clock (possibly a multiple of a GPS clock), and reports its every consumption over the $\tau$'th clock 
period as $D_i(\tau_i)$ to a centralised authority. This centralised authority aggregates this energy consumption and broadcasts a signal to the vehicles depending 
on whether the aggregated consumption exceeds $E(k,t)$ or not.  The probability whether the $i$'th vehicle travels in fully electric mode in the $\tau_i+1$'th period depends on this broadcasted signal. In the following two {\em use cases} we shall see examples of signals that can be used to orchestrate the fleet behaviour.

\subsection{Use case 1: Fair energy consumption}
The fair energy consumption case refers to the case when all the vehicles participating to the SPONGE program participate in the same manner; namely., have the same probability to travel in EV mode. In the SPONGE case illustrated in Section \ref{Sponge_case}, a simple proportional controller can be used:
\begin{equation}
\begin{array}{ll}
\mbox{\textbf{if}} & \sum_{i=1}^{N} D_i(k) < E_{av}(k+1)\\
\mbox{\textit{then}} & p_{i}^{EV}(k) = g_1(E_{av}(k+1) - \sum_{i=1}^{N} D_i(k)) 
\end{array}
\label{Exceed_SPONGE}
\end{equation}
in this case, at every interval of time (e.g., every minute), a vehicle travels in EV mode with a probability $p_{i}^{EV}(k)$ which is an increasing function ($g_1(\cdot)$) of the gap between the desired target of energy $E_{av}(k+1)$ and the currently available space in the vehicles $\sum_{i=1}^{N} D_i(k)$. Note that if enough space is already available, then the vehicles are allowed to travel in any way they desire. Also note that, as already anticipated, even if all $p_{i}^{EV}(k)$'s are set to 1, the goal might not be accomplished if not enough vehicles are travelling in the time interval of interest.\\
\newline
{\bf Comment 4: } Note that the SPONGE problem usually takes place on a day-scale. For instance, vehicles are scheduled to spend a given quantity of energy during the day, and are then recharged at night time, when idle, to refill the batteries. However, in a practical scenario, it is more convenient to match the energy over a number of time windows during the day. This has a number of benefits: if the match occurs after a few hours, then the cars travelling in the afternoon are excluded from the programme; on the other hand, if we split the matching problem in several windows of time, then every single car, travelling at any time, can be equally involved in the programme. Also, as soon as a new time window starts, then the matching problem can be adjusted taking into account new weather forecasts, if available, and whether the optimisation problem was feasible or not in the previous time window. Accordingly, $k$ refers to short time window, e.g., one minute, in Problem (\ref{Exceed_SPONGE})\\
\newline
{\bf Comment 5: } Note that the request that vehicles have to travel in EV mode with probability 0.6 can be implemented in practice either by making 60\% of the vehicles travel in EV mode, or by making 60\% of the traction provided by the EV engine, and the residual by the ICE engine, in every car.\\
\newline
As for the exact SPONGE case illustrated in Section \ref{Exact_Sponge_Case}, then a simple Proportional-Integral (PI) controller can be adopted in the following manner:
\begin{equation*}
\begin{array}{ll}
\mbox{\textbf{if}} & \sum_{i=1}^{N} D_i(k) < E_{av}(k+1)\\
\mbox{\textit{then}} & p_{i}^{EV}(k) = g_2(E_{av}(k+1) - \sum_{i=1}^{N} D_i(k))\\
\mbox{\textbf{elseif}} & \sum_{i=1}^{N} D_i(k) \geq E_{av}(k+1)\\
\mbox{\textit{then}} & p_{i}^{EV}(k) = 0, \forall i=1,...,N\\
\end{array}
\end{equation*}
the main difference in this case, is that the control objective is to \textit{exactly} deplete the batteries of the quantity $E_{av}(k+1)$, while vehicles are not allowed to over-deplete their batteries. Although such a solution may appear to penalise for the drivers (i.e., they are forced to travel in ICE mode to avoid over-depleting their batteries), it is very convenient for the grid, as it is possible to predict in advance exactly how much energy will have to be delivered to the fleet of vehicles. Also, as already mentioned in Comment 2, there may be good reasons for drivers to preserve a store of electric energy. \\
\newline
Note that the proposed approaches can be used to tackle many practical scenarios of interest. For instance, we can assume that a company provides a free battery-charging service to the PHEVs of the employees whenever there is enough power generated from some connected solar/wind plants in the surroundings of the company buildings. Then, there should always be some battery available for recharging whenever there is available energy from natural sources. Then, the employees collaborate to equally make space in their batteries to absorb the forthcoming energy. However note that, although such a scenario gives rise to fair solutions, still personal constraints of single employees are not taken into account. In this perspective, the scenario can be made more complicated as described in the following subsection.

\subsection{Use case 2: Utility optimisation}
The third scenario illustrated in Section \ref{Utility_Functions} is different from the previous two, since different probabilities should be computed for different users, taking personal constraints into account. For instance, the decision to travel in EV mode represents a potential cost to the owner that can be represented as an increasing function of $D_i$, i.e., of the energy spent to travel in EV mode. 
%An example of candidate function is depicted in Figure \ref{Likelihood_Figure}.
%\begin{figure}[ht]
%\centering
%\includegraphics[width=\columnwidth]{figs/Likelihood_Figure.pdf}
%\label{Likelihood_Figure}
%\caption{A function that represents the discomfort of PHEV's owners when travelling in EV mode. For instance, the discomfort can correspond to the likelihood $l_i$ of not having electrical power in the battery when needed.}
%\end{figure}
For example, the cost of travelling in EV mode can be mathematically formulated as the likelihood $l_i$ to not have electrical power when it is needed, e.g., because some areas might be accessible only in EV mode for pollution reasons (for example, the so-called \textit{umweltzonen} in Germany\footnote{\url{http://gis.uba.de/website/umweltzonen/umweltzonen.php}}.) Note that a similar discussion can be made in terms of discomfort of travelling in ICE mode.
This scenario allows the central infrastructure to explicitly take into account personal needs of PHEVs' owners and there are many ways to solve the mathematical problem that arises. One possible way is to obtain the solution by formulating the problems as regulation problems with constraints, and to use these constraints to solve optimisation problems as they arise. In the next section we provide one possible solution based on AIMD algorithms \cite{Twinlin, Fabian}.

\section{Simulations}
\label{Simulations}
\begin{figure}[ht]
	\centering
		\includegraphics[width=0.8\columnwidth]{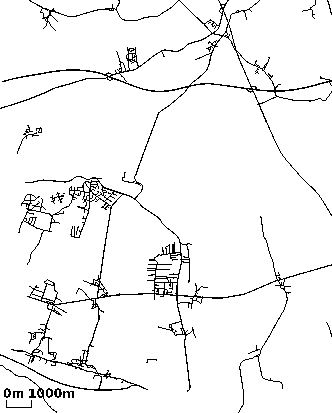}
	\caption{Road network in the Lower Saxony area in Germany used for our simulations, extracted from Open Street Map.}
	\label{fig:areaWithScale}
\end{figure}

\begin{figure*}[ht]
\begin{center}
\begin{tabular}[c]{ll}
\includegraphics[width=0.85\columnwidth]{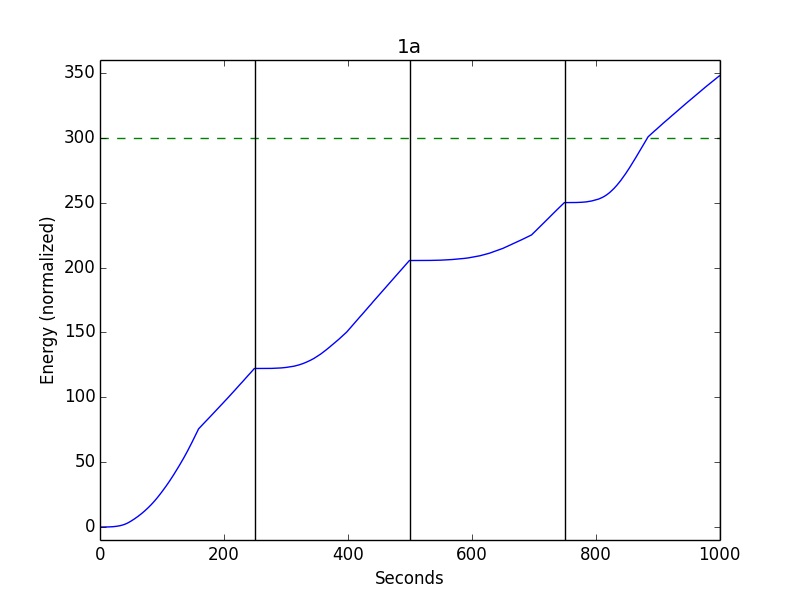} & 
\includegraphics[width=0.85\columnwidth]{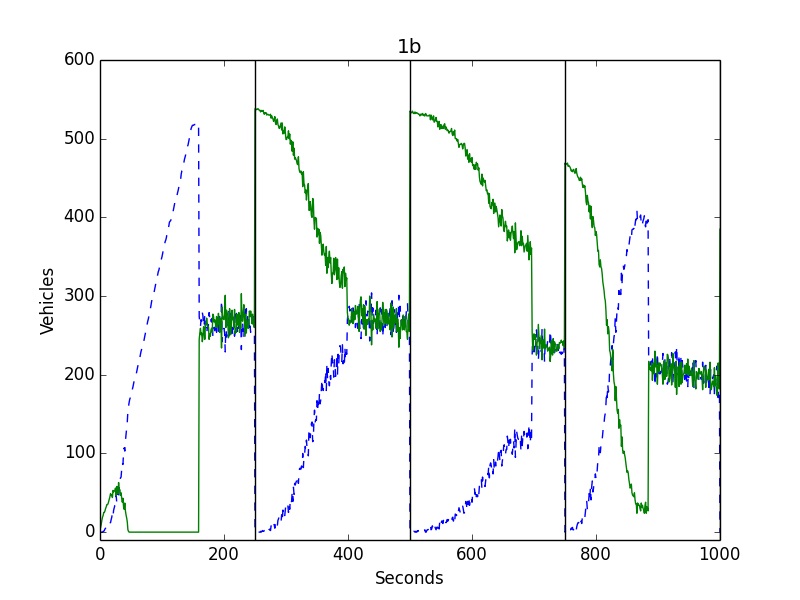} \\
\includegraphics[width=0.85\columnwidth]{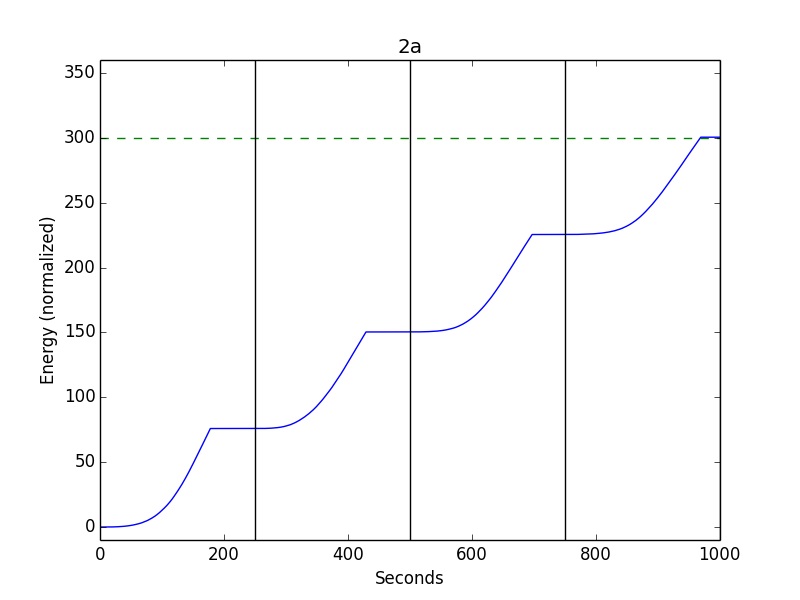} & 
\includegraphics[width=0.85\columnwidth]{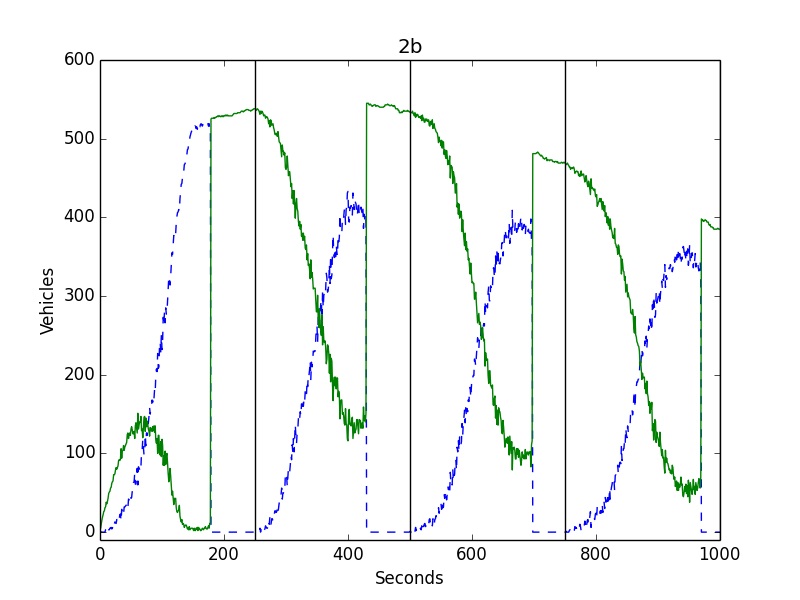} \\
\end{tabular}
\caption{Figures on the left show the total space for energy in the fleet of PHEVs in the two cases, in solid line, while the available energy from the grid is shown with the dashed line. Figures on the right show the mode in which the vehicles drive in the two cases to achieve the desired goal (dashed blue is EV mode, solid green is ICE mode). The vertical lines separate the time windows.}
\label{Simulation_Results}
\end{center}
\end{figure*}

We now present brief simulation results to show the efficacy of the proposed idea. The following simulations are performed using the popular mobility simulator SUMO~\cite{Sumo} and the given TRACI interface. A map of a rural area near Hamburg, Germany, was extracted from Open Street Map to be used as the underlying street network, and is shown in Figure~\ref{fig:areaWithScale}. Figure~\ref{Simulation_Results} shows the simulation results for the the first two algorithms. They refer to a time period of 1000 seconds (i.e., about 17 minutes). There are 4 time windows of 250 seconds each, and there are about 600 PHEVs on the road. In our simulations, we assume that vehicles that are running out of fuel, or whose battery is getting close to physical constraints (e.g., 10\% of the state of charge) are automatically discarded from the SPONGE programme. The simulation refers to a very simple example, and might correspond to the case when employees go to work using their PHEV vehicles, and the infrastructure regulates the driving mode in order to meet the target of energy that will be available at the workplace to recharge the vehicles.

\begin{figure*}[ht]
\begin{center}
\begin{tabular}[c]{ll}
\includegraphics[width=0.9\columnwidth]{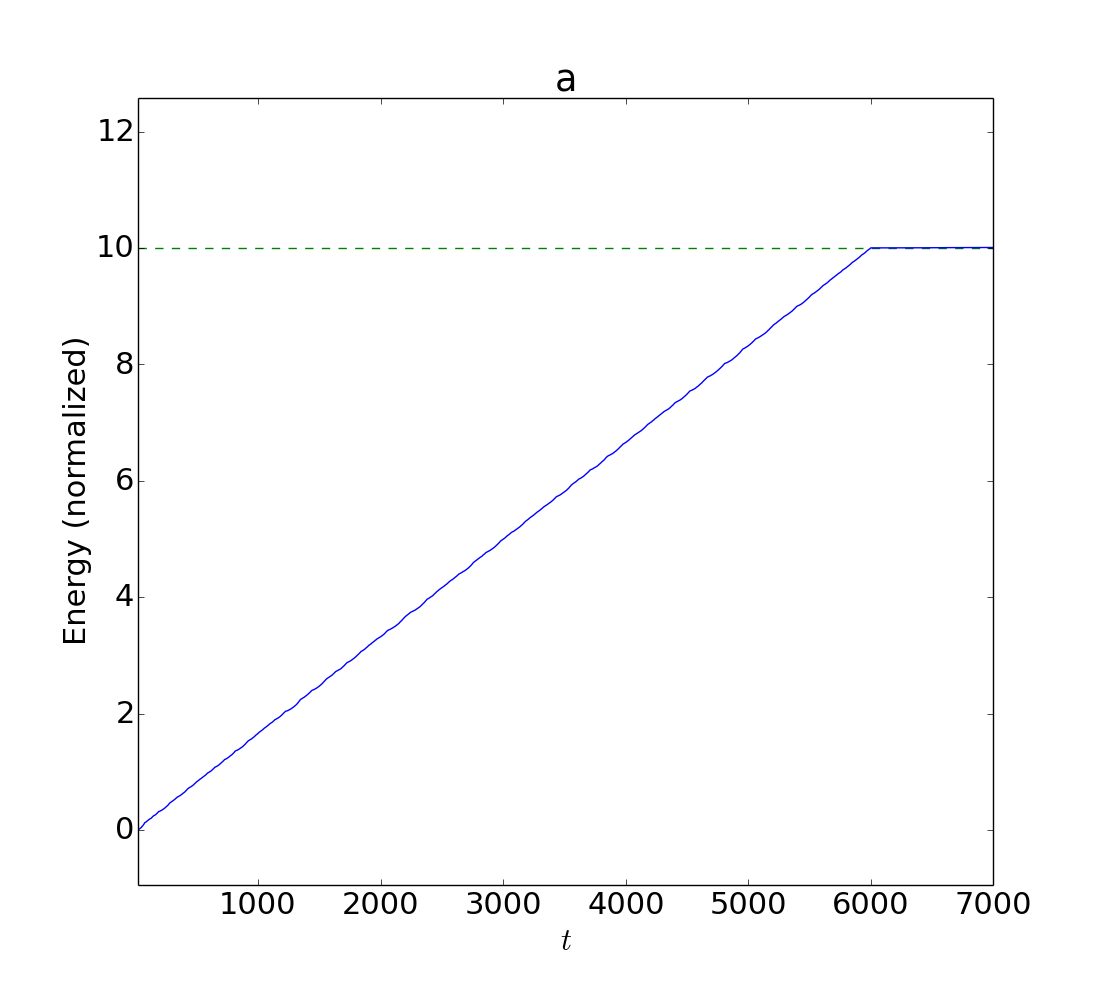} & 
\includegraphics[width=0.9\columnwidth]{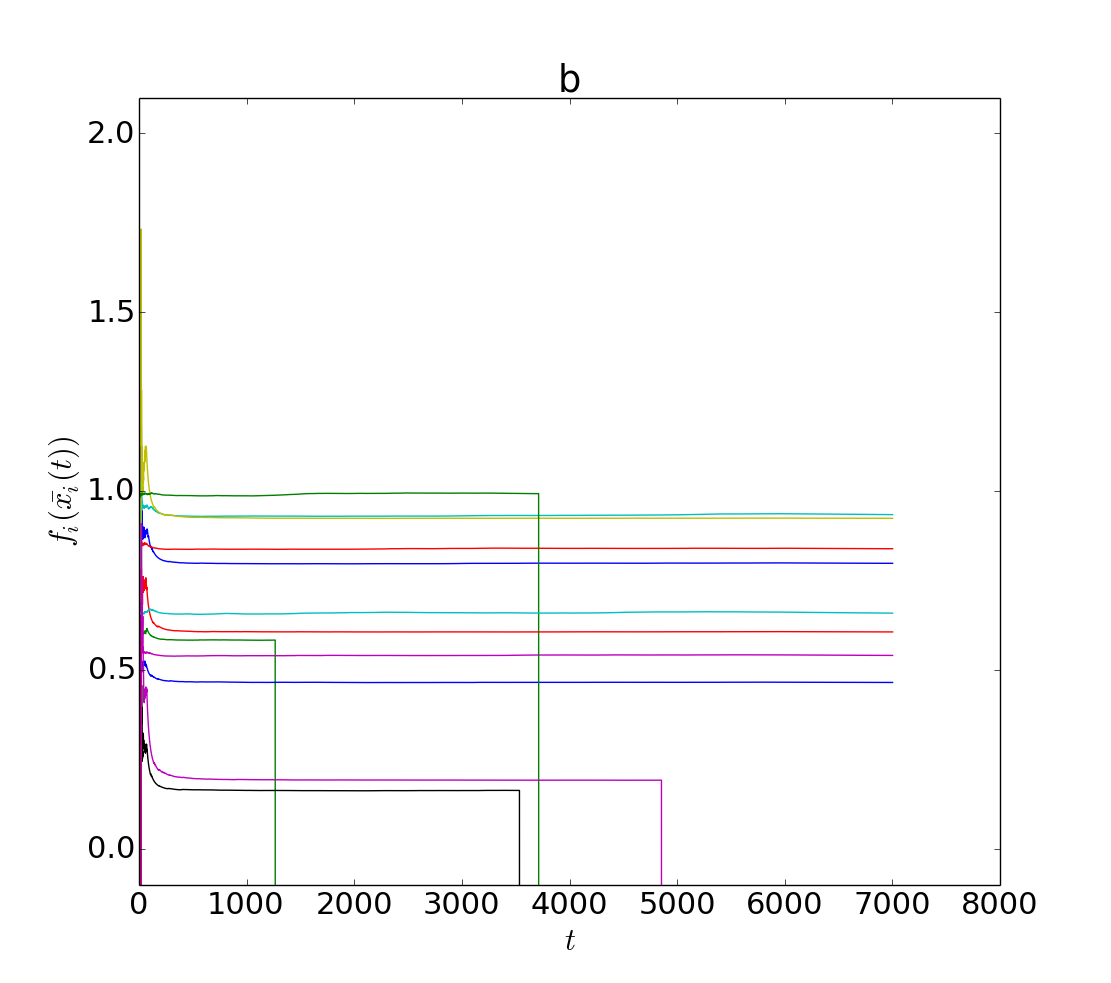} \\
\includegraphics[width=0.9\columnwidth]{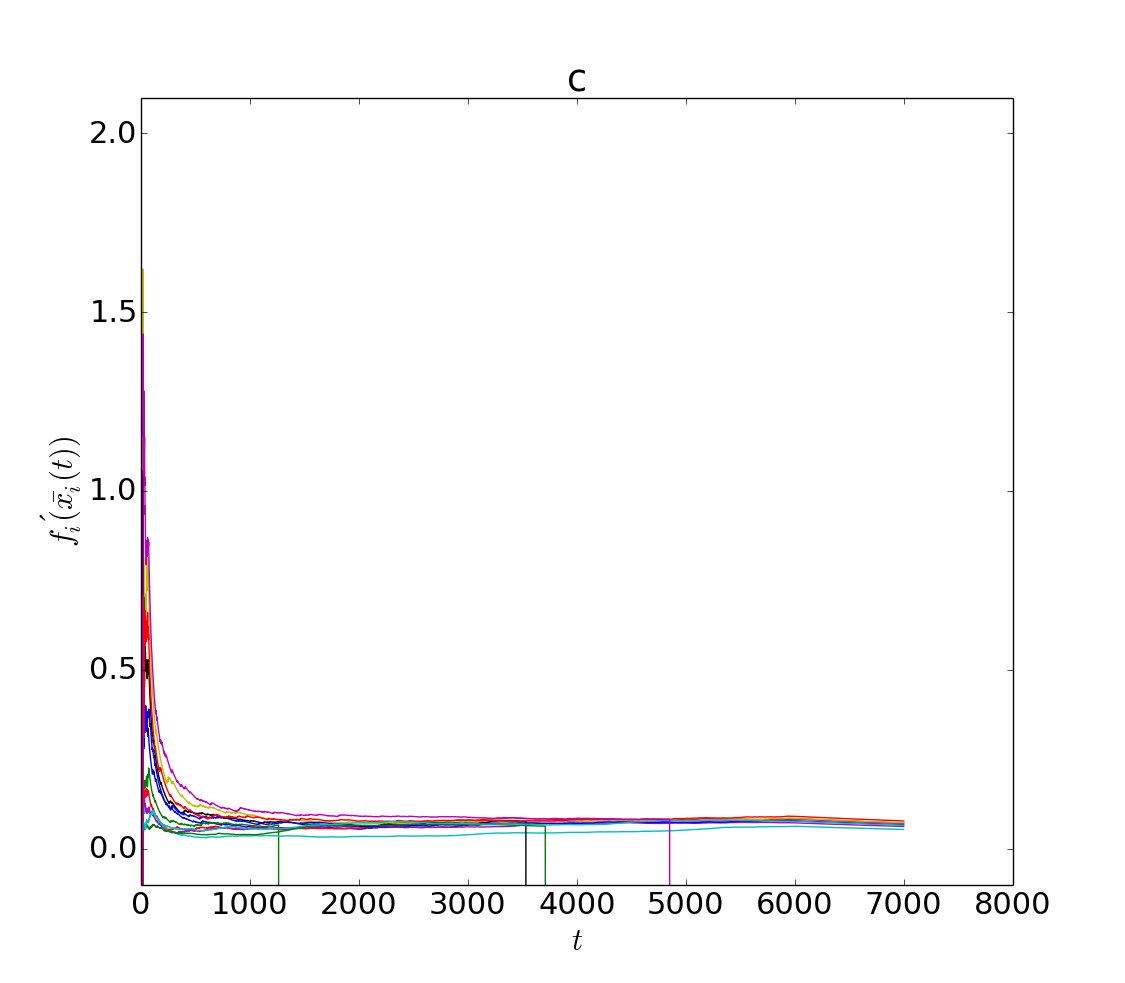} & 
\includegraphics[width=0.9\columnwidth]{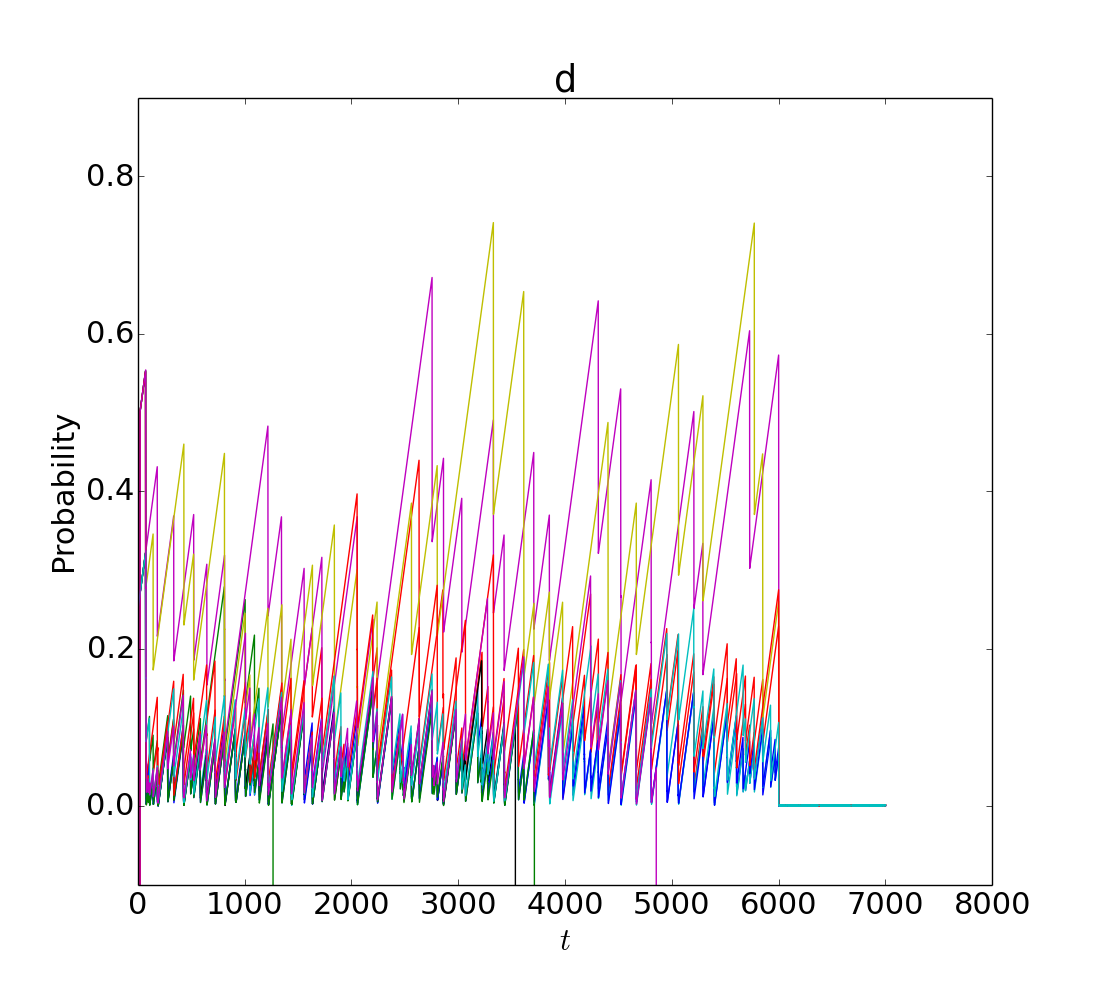} \\
\end{tabular}
\caption{Figure a shows that the overall constraint is again exactly satisfied at the end of the simulation. However, it is now achieved by minimising the sum of utility functions (single utility functions are shown in b). Figure c shows that the optimal solution has been achieved in fact, since there is a consensus on the derivative of the utility functions. Figure d emphasises that the utility optimisation problem is solved by assigning different probabilities of travelling in EV mode to different vehicles, according to their utility functions. Note that some vehicles reaching their destination before the end of the simulation exit the SPONGE programme in advance.}
\label{Simulation_Results AIMD}
\end{center}
\end{figure*}

We assume that vehicles travel in ICE mode when entering the simulation. In parallel, the infrastructure broadcasts simple indications that randomly make some vehicles travel in EV mode to free some space in the batteries. Once the target energy is matched, vehicles are allowed to travel in the mode they prefer (Figures~\ref{Simulation_Results}.1a and \ref{Simulation_Results}.1b) or in ICE mode in the exact SPONGE case (Figures~\ref{Simulation_Results}.2a and \ref{Simulation_Results}.2b). We made the assumption that when vehicles are allowed to choose their travelling mode, they choose with equal probability to travel in ICE or in battery mode. Note that the same pattern is repeated in every time window, until at the end of the simulation the overall target is met. The advantage of dividing the 1000 seconds into a number of smaller time windows of 250 seconds, is that both the vehicles that start their journey at the beginning and those at the end of the entire time frame, participate to the SPONGE programme.\\
\newline
The simulation results of the third scenario (utility optimisation) are shown in Figures~\ref{Simulation_Results AIMD}. We remind that in this scenario the exact equality between freed space and expected forthcoming energy (Figure~\ref{Simulation_Results AIMD}.a) is achieved by assigning different probabilities to travel in EV mode to different vehicles, according to some utility functions. We assumed that the inconvenience of vehicles in travelling in EV mode could be described through a convex quadratic function $f_i(\bar{x_i})=a_i \bar{x_i}$, with $\bar{x}_i$ being the share of time running in EV mode up to the current simulation step. Parameters $a_i$ were different for every vehicle, and in our simulation they were randomly chosen in the interval [0,1]. The evolution of the utility functions of some randomly selected vehicles is shown in Figure~\ref{Simulation_Results AIMD}.b. Note that the optimal solution of the Problem (\ref{Utility_Maximisation_Problem}) can be obtained by solving a consensus problem on the derivative of the utility functions (more detailed mathematical details together with a convergence proof can be found in \cite{Fabian}). Figure \ref{Simulation_Results AIMD}.c shows that the utility functions do indeed converge to the same value. Finally, Figure \ref{Simulation_Results AIMD}.d shows that the optimal solution is obtained by giving a different probability to travel in EV mode to each vehicle. The AIMD algorithm makes the probability of each vehicle to travel in EV mode linearly increase until the constraint is matched (congestion event). At that point, some probabilities decrease (back-off) in a multiplicative fashion to keep satisfying the constraint. Vehicles participates to the back-off step with a probability that is proportional to the derivative of their own utility function divided by the argument of the utility function (i.e., $\propto f'_i(\bar{x_i})/\bar{x_i}$). In this way, the optimal solution is obtained in a distributed way (i.e., without requiring vehicle-to-vehicle communication, or vehicle-to-infrastructure communication). The only communication requirement is from the infrastructure, that has to broadcast to all vehicles when the congestion event occurs. More details on the mathematical theory of the aforementioned algorithms, and more examples as well can be found in \cite{Fabian} and \cite{mingming}, and are here omitted due to space limits.

\section{Conclusions}
\label{Conclusions}

In this paper we have presented a new idea that takes advantage of the ability of PHEVs to both travel in electric and in fuel mode to absorb naturally generated electrical energy in a smart manner from the grid. From a theoretical perspective, such a problem can be easily formulated and solved using well-known algorithms for sharing a task among a number of distributed agents, e.g., AIMD algorithms \cite{mingming,Fabian}. From a practical point of view, note that the technology to remotely control the driving mode is also already available, as it was developed in \cite{Twinlin} for different purposes.\\
\newline Our current plan is to extend the preliminary simulation results given in Section \ref{Simulations} to more realistic and large-scale examples. In parallel, we intend to start implementing the approach in a reduced number of PHEVs, as a proof-of-concept of the paper idea. We shall adapt the experimental set-up of \cite{Twinlin} to the new case of interest, to remotely control the EV/ICE engine switching. Also, we shall integrate a reliable weather forecast software in the overall system, in order to take optimal decisions about when to switch from one mode to another mode.


\begin{thebibliography}{99}

\bibitem{Marinelli}
M.~Marinelli, F.~Sossan, G.T.~Costanzo and H.W.~Bindner, Testing of a Predictive Control Strategy for Balancing Renewable Sources in a Microgrid, {\it IEEE Transactions on Sustainable Energy}, vol. 5, no. 4, 2014.

\bibitem{Meibom}  
P.~Meibom, K.B.~Hilger, H.~Madsen and D.~Vinther, Energy Comes Together in Denmark: The Key to a Future Fossil-Free Danish Power System, {\it IEEE Power and Energy Magazine}, vol. 11, no. 5, pp. 46-55, 2013.

\bibitem{Bicik}  
V.~Bi\u{c}\'{i}k, O.~Holub, K.~Ma\u{r}\'{i}k, M.~Sikora, P.~Stluka and R.~D'hulst, Platform for coordination of energy generation and consumption in residential neighborhoods, {\it IEEE PES Innovative Smart Grid Technologies (ISGT) Europe}, 2012.

\bibitem{RAMP}
M.A.~Al Faruque, RAMP: Impact of rule based aggregator business model for residential microgrid of prosumers including distributed energy resources, {\it IEEE PES Innovative Smart Grid Technologies (ISGT)}, 2014.

\bibitem{Economist}
The~Economist, Grid-scale storage - Smooth operators, {\it The Economist Quarterly}, Q4, December 2014.

\bibitem{Liu}
C.~Liu, K.T.~Chau, D.~Wu and S.~Gao, Opportunities and Challenges of Vehicle-to-Home, Vehicle-to-Vehicle, and Vehicle-to-Grid Technologies, {\it Proceedings of the IEEE}, vol. 101, no. 11, pp. 2409-2427, 2013.

\bibitem{Tushar}
M.~Tushar, C.~Assi, M.~Maier and M.~Uddin, Smart Microgrids: Optimal Joint Scheduling for Electric Vehicles and Home Appliances, {\it IEEE Transactions on Smart Grid}, vol. 5, no. 1, pp. 239-250, 2014.

\bibitem{Griggs}
S.~St\"{u}dli, W.~Griggs, E.~Crisostomi and R.~Shorten, On Optimality Criteria for Reverse Charging of Electric Vehicles, {\it IEEE Transactions on Intelligent Transportation Systems}, vol. 15, no. 1, pp. 451-456, 2014.

\bibitem{PNAS}
C.W.~Tessum, J.D.~Hill and J.D.~Marshall, Life cycle air quality impacts of conventional and alternative light-duty transportation in the United States, {\it PNAS} 2014.

\bibitem{CNR_IEVC}
E.~Ancillotti, R.~Bruno, E.~Crisostomi and M.~Tucci, Using Electric Vehicles to Improve Building Energy Sustainability, {\it IEEE International Electric Vehicle Conference (IEVC)}, 2014.

\bibitem{Edash}
e-DASH, Electricity Demand and Supply Harmonisation for Electric Vehicles, available online at \url{http://edash.eu/} [last checked, December 2014].

\bibitem{Twinlin}
A.~Schlote, F.~H\"{a}usler, T.~Hecker, A.~Bergmann, E.~Crisostomi, I.~Radusch and R.~Shorten, Cooperative Regulation and Trading of Emissions Using Plug-in Hybrid Vehicles, {\it IEEE Transactions on Intelligent Transportation Systems}, vol. 14, no. 4, pp. 1572-1585, 2013.

\bibitem{Wind}
B.-M.~Hodge, D.~Lew, M.~Milligan, H.~Holttinen, S.~Sillanp\"{a}\"{a}, E.~G\'{o}mez L\'{a}zaro, R.~Scharff, L.~S\"{o}der, X.G.~Lars\'{e}n, G.~Giebel, D.~Flynn and J.~Dobschinski, Wind Power Forecasting Error Distributions: An International Comparison, {\it 11th Annual International Workshop on Large-Scale Integration of Wind
Power into Power Systems as well as on Transmission Networks for Offshore Wind Power Plants}, Lisbon, Portugal, 2012.

\bibitem{Sun}
J.~Zhang, B.-M.~Hodge, A.~Florita, S.~Lu, H.F.~Hamannk and V.~Banunarayanan, Metrics for Evaluating the Accuracy of Solar Power Forecasting, {\it 3rd International Workshop on Integration of Solar Power into Power Systems}, London, England, 2013.

\bibitem{Fabian}
F.~Wirth, S.~St\"{u}dli, J.Y.~Yu, M.~Corless and R.~Shorten, Nonhomogeneous Place-Dependent Markov Chains, Unsynchronised AIMD, and Network Utility Maximization, available online at \url{http://arxiv.org/pdf/1404.5064v2.pdf}.

\bibitem{Oak_Ridge}
M.L.~Shaltout, D.~Chen, A.A.~Malikopoulos and S.~Pannala, Multi-Disciplinary Decision Making and Optimization for Hybrid Electric Propulsion Systems, {\it IEEE International Electric Vehicle Conference (IEVC)}, Florence, Italy, 2014.

\bibitem{deHoog}
J.~de Hoog, T.~Alpcan, M.~Brazil, D.A.~Thomas and I.~Mareels, A Market Mechanism for Electric Vehicle Charging Under Network Constraints, {\it IEEE Transactions on Smart Grid}, accepted for publication, 2014.

\bibitem{Sonja}
S.~St\"{u}dli, E.~Crisostomi, R.~Middleton and R.~Shorten, Optimal real-time distributed V2G and G2V management of electric vehicles, {\it International Journal of Control}, vol. 87, no. 6, pp. 1153-1162, 2014.

\bibitem{mingming}
E.~Crisostomi, M.~Liu, M.~Raugi and R.~Shorten, Plug-and-play Distributed Algorithms for Optimised Power Generation in a Microgrid, {\it IEEE Transactions on Smart Grid}, vol. 5, no. 4, pp. 2145-2154, 2014.

\bibitem{Sumo}
D. Krajzewicz, M. Bonert, P. Wagner, The open source traffic simulation package SUMO, {\it RoboCup 2006 Infrastructure Simulation Competition}, Bremen, Germany, 2006.

\end{thebibliography}
\end{document}